\documentclass[12]{amsart}
\usepackage{bbm}
\usepackage{fullpage}
\title{On standard Young tableaux of bounded height}
\author{M J Mishna}

\usepackage{amsmath, amsthm, amsfonts}
\usepackage{verbatim}
\usepackage{caption}
\usepackage{subcaption}
\usepackage{pifont, xcolor, tikz}
\usepackage{graphicx}
\usepackage{ytableau}
\usepackage{array}
\usepackage{amssymb}
\usepackage{enumitem}

\newcommand{\yt}[1]{
\ensuremath{\begin{array}{c}
\begin{ytableau}#1\end{ytableau}
\end{array}}
}

\newtheorem{theorem}{Theorem}
\newtheorem{proposition}[theorem]{Proposition}

\newtheorem{conjecture}[theorem]{Conjecture}

\newtheorem{corollary}[theorem]{Corollary}

\theoremstyle{remark}

\newtheorem{example}{Example}

\newcommand{\bz}{\ensuremath{\mathbf{z}}}

\usepackage{pifont, xcolor, tikz}

\usetikzlibrary{positioning}
\usetikzlibrary{decorations.pathreplacing}

\tikzset{partition/.style={fill,circle,inner sep=1pt},
         part/.style={baseline=0,scale=0.5,bend left=45},
         partlabel/.style={below}}
\usetikzlibrary{shapes,arrows}
\tikzstyle{pnt}=[draw,ellipse,fill,inner sep=1pt]

\tikzstyle{opnt}=[ ]
\tikzstyle{pntt}=[draw,ellipse,fill,inner sep=0.5pt]
\tikzstyle{point}=[draw,ellipse,fill,inner sep=2pt]

\makeatletter
\newcommand{\oset}[2]{%
  {\mathop{#2}\limits^{\vbox to \ex@{\kern-\tw@\ex@
   \hbox{\ensuremath{#1}}\vss}}}}
\makeatother

\begin{document}

\begin{abstract}
  We survey some recent works on standard Young tableaux of
  bounded height. We focus on consequences resulting from numerous
  bijections to lattice walks in Weyl chambers. 
\end{abstract}
\maketitle

\section{Introduction}
\label{sec:Intro}
Standard Young tableaux are a classic object of mathematics, appearing
in problems from representation theory to bijective
combinatorics. Lattice walks restricted to cones are similarly a
fundamental family, and they encode a wide variety of combinatorial
structures from formal languages to queues. Standard Young
tableaux of bounded height are in bijection with several different
straightforward classes of lattice walks. This connection not only
elucidates several sources of ubiquity on both accounts, but
facilitates exact and asymptotic enumeration, as well as parameter
analysis. This survey describes the cross developments over the past
the past 30 years, and highlights some open problems. For more
background on the tableaux, we recommend the surveys~\cite{Saga90,
  AdRo15} and particularly the article of Stanley~\cite{Stan07}.

We begin by fixing our notation and conventions. Let
$\lambda=(\lambda_1, \lambda_2, \dots, \lambda_t)$ be a partition of
$n$ into $t$ parts. We write $\lambda\vdash n$ and $\ell(\lambda)=t$
in this case.  The \emph{Ferrers diagram of shape $\lambda$} is a
representation of $\lambda$ comprised of boxes indexed by pairs
$\{(i,j): 1\leq i \leq t; 1\leq j\leq \lambda_i\}$. Such a diagram is
of size~$n$, and height~$t$. A key parameter in our study is the
number of columns of odd length.

A \emph{standard Young tableau\/} of size $n$ is a filling of a
Ferrers diagram using precisely the integers $1$ to $n$. The entries
strictly increase to the right along each row and strictly increase
down each column.  A tableau is \emph{semi-standard} if the entries
weakly increase along each row and strictly increase down each
column. In this work we are interested in the number of standard Young
tableaux of a given size with the height is bounded by a fixed value.

\subsection{Enumeration formulas}
\label{sec:History}
Rather classically, the number of standard Young tableaux of shape
$\lambda$ is denoted~$f^\lambda$ and is given by the hook length formula:
\begin{equation}
\label{eqn:flambda}
f^\lambda = \frac{n!}{\prod_c h_c}\quad \text{where} \quad
h_c=\lambda_i+\operatorname{card}\{j:\lambda_j\geq i\} -i-j+1.
\end{equation}
The following formulation is due to MacMahon:
\[
f^\lambda = (\lambda_1+\dots+\lambda_d) ! \det\left(
  \frac{1}{(\lambda_i-i+j)!} \right)_{1\leq i,j\leq d}.
\]
The number of standard Young tableaux of height at most $k$ is thus the sum
\begin{equation}
\label{eqn:ykn}
y_k(n)\equiv\sum_{\substack{\lambda\vdash n\\
  \ell(\lambda)\leq k}} f^\lambda.
\end{equation}

\begin{figure}
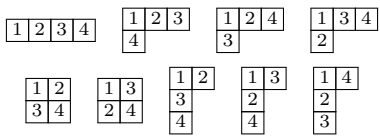

\ytableausetup{smalltableaux}
\begin{equation*}
\ytableausetup{centertableaux}
\begin{ytableau}
 1 & 2 & 3 & 4
\end{ytableau}\quad
\begin{ytableau}
 1 & 2 & 3\\
4
\end{ytableau}\quad
\begin{ytableau}
 1 & 2 & 4\\
3
\end{ytableau}\quad
\begin{ytableau}
 1 & 3 & 4\\
2
\end{ytableau}
\end{equation*}
\begin{equation*}
\begin{ytableau}
 1 & 2 \\
3 & 4
\end{ytableau}\quad
\begin{ytableau}
 1 & 3\\
2&4
\end{ytableau}\quad
\begin{ytableau}
1 & 2\\
3\\
4
\end{ytableau}\quad
\begin{ytableau}
 1 & 3\\
2\\
4
\end{ytableau}\quad
\begin{ytableau}
 1 & 4\\
2\\
3
\end{ytableau}
\end{equation*}
\caption{All standard Young tableaux of size 4 and height at most 3}
\label{fig:y34is9}
\end{figure}

The first enumerative formulas for Young tableaux where the height is
an explicit consideration appear in the 1960s, when Gordon and Houten
studied $k$-rowed plane partitions whose non-zero parts strictly
decrease along rows and columns, in addition to some related
variants. In their series of Notes on Plane
Partitions~\cite{GoHo68,GoHo69}, they give some formulas for the
generating functions in terms of infinite products and determinants.
Regev~\cite{Rege81} first determined exact expressions for $y_2(n)$
and $y_3(n)$. The formulas are
equivalent to the following. Here~$C_k$ denotes the $k$-th
\emph{Catalan number\/}\footnote{$C_n\equiv\binom{2n}{n}\frac{1}{n+1}$}:
\begin{equation}\label{eqn:y3n}
y_2(n) = \binom{n}{\lfloor n /2\rfloor}\qquad
y_3(n)=\sum_{k=0}^{\lfloor n /2\rfloor} \binom{n}{2k}C_k.
\end{equation}
The numbers $y_3(n)$ are also known as \emph{Motzkin
  numbers}. Figure~\ref{fig:y34is9} illustrates the nine standard
Young tableaux of size four and of height at most three.  Around the
same time, Gessel~\cite{Gess90} found an expression for~$y_4(n)$, and
Gouyou-Beauchamps~\cite{Gouy89} found the following expressions for
$y_4(n)$ and $y_5(n)$:
\begin{equation}\label{eqn:y4n}
\begin{split}
y_4(2n)=C_nC_n \qquad y_4(2n+1)=C_nC_{n+1}\\
y_5(n)=\sum_{i=0}^{\lfloor n/2 \rfloor}
\frac{3!n!(2i+2)!}{(n-2i)!i!(i+1)!(i+2)!(i+3)!}.
\end{split}
\end{equation}
No comparable expression for $y_6(n)$ has appeared in the literature.
The presence of binomials in general, and Catalan numbers in
particular, is a strong hint that these tableaux are related to well
understood combinatorial classes. 

\subsection{The exponential generating function}
We study $y_k(n)$ via $Y_{k}(t)$, the exponential generating function
for $y_k(n)$:
\begin{equation}
Y_{k}(t) \equiv \sum y_k(n) \frac{t^n}{n!}.
\end{equation}
The formulas depends on the parity of the height. The formula for $Y_{2k}(t)$ was obtained by Gordon~\cite{Gord71} by
reducing a Pfaffian of Gordon and Houten~\cite{GoHo68}, and
Gessel~\cite{Gess90} found the formula for odd heights. They are both expressed in
terms of the \emph{hyperbolic Bessel function of the first kind of
  order~$j$}\[b_j\equiv
  I_j(2t)=\sum_{n=0}^\infty\frac{t^{2n+j}}{n!(n+j)!}.\] The
formulas are:
\begin{align}
\label{eqn:Ykt}
Y_{2k}(t)&= \det\left(b_{i-j}+b_{i+j-1}\right)_{1\leq i,j\leq k}\\
Y_{2k+1}(t)&=e^t\ \!\det\big(b_{i-j}-b_{i+j}\big)_{1\leq i,j\leq k}. 
\end{align}

For example, $Y_2(t)= b_0+b_1$ and
$Y_4(t)=b_0^2+b_0b_1+b_0b_3-2b_1b_2-b_2^2-b_1^2+b_1b_3.$ These
expressions grow very fast as polynomials in $b_j$, however, they are
amenable to some further study including asymptotic analysis of the
coefficients.  They also imply that $y_k(n)$ can be expressed by
binomial sums, although an expression that one could compute may be
too complex to be of any use.

One important consequence of these formulas is that they resolved a
question that Stanley~\cite{Stan80} had asked almost 10 years before
these formulas appeared about the nature of the function
$Y_k(t)$. Specifically, he asked if $Y_k(t)$ is \emph{D-finite}, that
is, does it satisfy a differential equation with polynomial
coefficients. Gouyou-Beauchamps showed that $Y_4(t)$ is D-finite from
his formula, and Gessel~\cite{Gess90} proved this result for
general~$k$ since the expressions are a polynomial
combination of Bessel functions (possibly times an exponential).
Bessel functions are D-finite so the result follows from closure
properties. 
\subsection{Schur Functions}
To prove these formulas Gessel started with another important formula
for standard Young tableaux.  Schur functions can be described via a summation over the set of
all semi-standard Young tableaux of shape $\lambda$:
\begin{equation}
s_\lambda = \sum_{T\in \operatorname{SSYT}(\lambda)}
x^{\operatorname{content}(T)} = \sum_{T\in \operatorname{SSYT}(\lambda)}x_1^{t_1}x_2^{t_2}\dots x_k^{t_k}.
\end{equation}
The exponents describe the content of the tableau: $t_i$ is the number
of occurrences of $i$ in $T$.  Thus, the coefficient of the monomial
$x_1\dots x_n$ in this expression is the number of standard Young
tableaux of shape $\lambda$.  We deduce the formula
\[
y_k(n)= \sum_{\substack{\lambda\vdash n\\ \ell(\lambda)\leq k}} [x_1x_2\dots x_n] s_\lambda 
=[x_1x_2\dots x_n]\sum_{\substack{\lambda\vdash n\\ \ell(\lambda)\leq k}} s_\lambda .
\]

This kind of coefficient extraction in symmetric functions can be
framed as a homomorphism. This was done by Gessel in his PhD thesis
(Theorem 3.5) and also by Jackson and Goulden~\cite[Lemma
4.2.5]{GoJa83}.  In the case of the homogeneous complete symmetric
function $h_n$, it is easy to see that
$[x_1x_2\dots x_n]h_k =\mathbbm{1}_{n=k}$. One would like to apply
this to the \emph{Jacobi-Trudi\/} identity, which is an expression for
a Schur function in terms the homogeneous complete symmetric
functions:
\[s_\lambda=\det\left(h_{\lambda_i+j-i}\right)^{\ell(\lambda)}_{i,j=1}.\] 

The truth is slightly more complicated. Indeed Gordon and Houten did
much of the heavy lifting for this problem and had expressed the
number of semi-standard tableaux as a determinant of homogenous
complete symmetric functions. Gessel extracted the generating function
for standard Young tableaux, and derived the formulas in
Equation~\eqref{eqn:Ykt}.

\subsection{The Robinson-Schensted correspondence}
The following identity is incredibly evocative to combinatorialists:
\begin{equation}
\sum_{\lambda\vdash n} (f^{\lambda})^2 =n!
\end{equation}
The bijective correspondence between pairs of
standard Young tableaux of the same shape and permutations was
described by Robinson in the 1930s and also by Schensted in the
1960s. Below we describe the Schensted algorithm which builds a
pair of tableaux by parsing the permutation and
incrementally building two tableaux. We refer readers to
Sagan's book~\cite{Saga01} for additional details. 

An algorithmic description of the bijection begins with a permutation
$\sigma \in \mathbb{S}_n$ and a pair of empty tableaux, denoted
$(P_0,Q_0)$. For $i$ from 1 to $n$ we create $P_i$ by adding a box
with entry $\sigma(i)$ to $P_{i-1}$ via a special insertion
algorithm. For each $i$ we add precisely one box. The position is
noted and a box with entry $i$ in the same location is added to
$Q_{i-1}$. The sequence of $Q$ tableaux record the history of the box
additions. The pair $(P_n, Q_n)$ is returned.

The row insertion takes as input a possibly incomplete standard Young
tableau and adds an integer $m$ not already in the tableau. The
process acts incrementally along each row. If $m$ is bigger than all
of the elements in the row under consideration, it is placed at the end
of it. If not, it finds its natural place and bumps the larger
value. The bumped value is inserted in the tableau formed by the lower
rows by the same process.

\begin{example}
Let $\sigma$ be the involution $\sigma=
7\,2\,9\,6\,10\,4\,1\,8\,3\,5$. 
If the Schensted algorithm is applied, the penultimate step
returns the following two tableaux:
\[P_9=\begin{ytableau}
1 & 3 & 8\\
2 &  4& 10\\
6 & 9 \\
7
\end{ytableau}\quad
Q_9=\begin{ytableau}
1 & 3 & 5\\
2 &  4& 8\\
6 & 9 \\
7
\end{ytableau}.\] The final step is to insert 5 into $P_9$. The 5 fits
between the 3 and 8 in the first row, hence it bumps 8 to the next row
which then bumps $10$ which settles at the end of row three. The
position of the new square is $(3,3)$, and the algorithm finishes by
adding this square to $Q_9$, with an entry value of 10.

The end result is 
\[P_{10}=Q_{10}=
\ytableausetup{smalltableaux}
\begin{ytableau}
1 & 3 & 5\\
2 &  4& 8\\
6 & 9 &10\\
7
\end{ytableau}\]
The involution has two  fixed points (2 and 8)  and the resulting tableaux each have
two columns of odd length (specifically, length 3).
\end{example}

The example illustrates some properties which are true in general. The
length of the longest increasing subsequence of~$\sigma$ is equal to
the length of the first row of the tableaux. If $\sigma$ is an
involution, then the number of fixed points of $\sigma$ equals the
number of columns of odd length in $\lambda$. That standard Young
tableaux of size~$n$ are equinumerous with involutions was known to
Frobenius and Schur.

Viennot~\cite{Vien78} described a very beautiful geometric
construction using \emph{shadow lines} to give a more intuitive
illustration of this slightly mechanical bijection. It should be
better incorporated into the bijections we encounter in later
sections. A second geometric construction is given by Fomin using
growth diagrams~\cite{fomi95}. The construction can be mined for additional
information.

\subsection{Plan of the article}
The first formulas found $Y_k(t)$ were recognized to resemble
generating functions for walks in Weyl chambers deduced by Grabiner
and Magyar~\cite[Section 6.2]{GrMa93}. We examine
the relevant developments made in the early 1990s~\cite{Wilf92,
  Wilf95, GeWeWi98, Zeil07, Xin10} in Section~\ref{sec:LatticePaths}.

Bijective proofs of some of these connections are more
recent~\cite{Euetal13, Buetal16, Krat16}, although the work of
Gouyou-Beauchamps dates back to the late 1980s. All of these authors'
proofs pass through secondary objects, such as coloured Motzkin paths,
or matchings. We consider these in Section~\ref{sec:Collection}.

One of the most common strategies for enumerating lattice walks
restricted to cones with symmetry in the set of allowable steps
involves a sub-series extraction from a rational function.
Remarkably, many combinatorial classes with transcendental D-finite
generating functions share this property. It is an open problem to
answer under which conditions this might be universally true. In part,
it is a useful formulation since in ideal cases we can answer
questions about the order and complexity of the recurrences, and also
(re-) derive asymptotic formulas. We describe such extractions, and
the resulting tableaux generating functions in
Section~\ref{sec:diag}.

From the Robinson-Schensted correspondence we see that the
combinatorial cousin to the standard Young tableau of bounded height
is the permutation with bounded longest increasing subsequence. Many
of the techniques described here can also enumerate these classes. For
example, Gessel determined expressions for the generating functions,
and they also have lattice walk interpretations. We summarize
results in Section~\ref{sec:permutations}.

We conclude with some natural generalizations and open problems.

\begin{figure}\center
\ytableausetup{smalltableaux}
$\iota=\begin{ytableau}
 1 & 3 & 5\\
2 &  4& 8\\
6 & 9 &10\\
7
\end{ytableau}$\\[2em]

$\sigma=\begin{array}{llllllllll}
1&2&3&4&5&6&7&8&9&10\\
7&{\color{red}2}&9&6&10&{\color{red}4}&1&{\color{red}8}&3&5\\
\end{array}$

 \ytableausetup{smalltableaux}

$$\omega=(b_1, b_2, b_1, b_2, b_1, b_3, b_4, b_2, b_3, b_3)$$
\smallskip

$\mu=$\begin{tikzpicture}[scale=0.4]
\foreach \x in {1,2,...,10}{
\node[pnt, label=below:{$\x$}](\x) at (\x,0){};}
\draw[bend left=45](1) to (10);
\draw[bend left=45](3) to (6);
\draw[bend left=45](4) to (9);
\draw[bend left=45,color=gray](2) to (11,1);
\draw[bend left=45,color=gray](8) to (11,.5);
\draw[bend left=45](5) to (7);
\end{tikzpicture}
\smallskip

$$\theta={\tiny \left(\varnothing,  
\yt{\ },  \yt{\ \\ \ }, \yt{\  & \  \\ \ }, \yt{\  & \  \\ \ & }, 
\yt{\  & \ & \ \\ \ & \ }, \yt{\  & \ & \ \\ \ }, \yt{\ & \ \\ \ },
\yt{\  & \  \\ \ & }, \yt{\ & \ \\ \ },
 \yt{\ & \  } \right)}$$

\caption{Representatives of classes of objects in bijection with
  standard Young tableaux of bounded height, in particular their image
  of $\iota$, a standard
  Young tableaux of height 3 with two odd columns; $\sigma$ is an
  involution with maximal increasing sequence of length 3; $\omega$ is
a ballot walk; $\mu$ is an open arc diagram for a partial matching;
$\theta$ is an oscillating tableaux ending on a row shape.}
\label{fig:Collection}
\end{figure}

\section{Lattice Walk Models}
\label{sec:LatticePaths}
There are no fewer than five different lattice walk classes that are
equinumerous with standard Young tableaux of bounded height.  Most of
them are defined using $d$-dimensional Weyl chambers\footnote{For
  convenience we define the chambers using non-strict inequalities,
  our bijective statements can equivalently be given under strict
  inequalities, upon applying the coordinate shift
  $\widetilde{x}_i=x_i+d+1-i$.} of type~$C$, denoted by~$W_C(d)$
\[W_C(d)\equiv\left\{ (x_1,x_2,\dots,x_d): x_1 \geq x_2 \geq \dots
  \geq x_{d} \geq 0\right\}.\]

Let $\{e_1,\dots,e_d\}$ denote the standard basis of
$\mathbb{R}^{d}$. A lattice walk model is defined by a set of
allowable steps and a region which confines each walk. A walk is a sequence of steps.
For example, the set 
\[\mathcal{S}=\{\pm e_i: 1\leq i \leq d\} \]
defines the~$d$ dimensional \emph{simple}\/ step set. This is
considered in several different bounding cones. 

Gessel and Zeilberger~\cite{GeZe92} considered general walks in Weyl
chambers, and demonstrated an enumeration strategy for models where
the stepset possesses a certain kind of symmetry, and avoid jumping
over boundaries. Such walks have come be called walks
\emph{reflectable}. The generating function for reflectable walk
models with specified starting and endpoints (excursions) are written
as a coefficient extraction from a signed sum of unrestricted
walks. This idea, which appears frequently in lattice walk
enumeration, is a version of the reflection principle. The expressions
were made explicit for $W_C(d)$ by Grabiner and Magyar for the simple
step set, which happens to be reflectable.
\begin{theorem}[Grabiner and Magyar 1993~\cite{GrMa93}] For
  fixed~${\lambda}, {\mu} \in W_C(d)$, the exponential generating
  function $O_{{\lambda}, {\mu}}(t)$ of the simple walks in $W_C(d)$
  from $\lambda$ to $\mu$, counted by their lengths, satisfies
\[
O_{{\lambda}, {\mu}}(t) 
= \det \left( b_{{\mu_i}-{\lambda_j}}-b_{{\mu_i}+{\lambda_j}}\right)_{1\leq i, j\leq d}.
\]
\end{theorem}
This is immediately reminiscent of the generating function formulas in
Equation~\eqref{eqn:Ykt}. Several authors have made the connection, in
particular, Gessel, Weinstein and Wilf~\cite{GeWeWi98},
Zeilberger~\cite{Zeil07}, Xin~\cite{Xin10}, Eu~\emph{et
  al.}~\cite{Eu10, Euetal13}, Burrill~\emph{et al.}~\cite{Buetal16}
and Courtiel~\emph{et al.}~\cite{Coetal18}. In almost every case there
is a natural parameter which is equidistributed with number of odd
columns in the tableaux. We describe these classes next.

\subsection{Ballot walks}
We can build an integer sequence from the entries of a standard Young
tableaux. We define the associated lattice word
$w=(w_1, w_2, \dots, w_n)$ by setting $w_i=j$ if the entry $i$ is in
the $j$-th row.  This word has the property that for any prefix the
number of occurrences of $\ell$ is greater than or equal to the number
of occurrences of $\ell+1$ since the columns of the associated tableau
are strictly decreasing. These are also known as generalized Ballot
sequences.

We associate a step naturally to each letter:
\[ w_1\rightarrow e_1,\qquad w_i \rightarrow e_{i}-e_{i-1} \quad
  \text{for } (2\leq
  i\leq k-1), \qquad w_k \rightarrow -e_{k-1}. 
\]
Restricting walks to the first orthant is equivalent to the ballot
condition. That is, the ballot word of a tableaux of height at most
$k$ gives a natural encoding as a lattice walk in the cone
$\mathbb{R}^{k-1}_{\geq 0}$ with steps from the following stepset:
\[\mathcal{B}\equiv\{e_1, -e_{k-1}\}\cup \{e_{i}-e_{i-1}: 2\leq
  i\leq k-1\}. \]
Figure~\ref{fig:Collection} has an example.
\subsection{Lazy walks}
\label{sec:Lazy}
To define this step set, we set $\overline{0}$ to be the zero step
(whence the label ``lazy'').  Zeilberger~\cite{Zeil07} noted,
(although, he attributes the proof to Xin without citation) that the
number of excursions in $W_C(k)$ starting and ending at the origin
using the stepset $\mathcal{L}$, given by
\[
\mathcal{L}\equiv\{e_i, -e_i: 1\leq i\leq k \} \cup \{\overline{0}\}
\]
is $y_{2k+1}(n)$. A small computation suggests that the distribution of
the zero steps matches the distribution of odd columns in the Young
tableaux.  He remarks that it would be interesting to find a (bijective)
proof of this result. 

\subsection{Generalized Motzkin Paths}
\label{sec:Motzkin}
In Zeilberger's lazy walks, the $k=1$ case encodes the classic Motzkin
walks, consistent with the longstanding observation that $y_3(n)$ is
the number of Motzkin words. One could view the higher dimension lazy
walks as a generalization of Motzkin words, and this notion was first
formalized by Eu~\cite{Eu10}, and subsequently by Eu, Fu, Hou, and
Hsu~\cite{Euetal13}. They add a counter component, and describe an
explicit bijection between the Motzkin paths of length~$n$ and the
standard Young tableaux of size~$n$ with at most three rows.

To prove this, they considered the lazy lattice walks, and
then algorithmically mapped the steps in. The
odd and even cases reconcile as follows. 
\begin{theorem}[Eu, Fu, Hou, and Hsu 2013{~\cite[Theorem 1.1]{Euetal13}}]
Consider the lattice model defined by the step set 
\[\mathcal{M}\equiv\{e_1\}\cup\{e_1 +e_2,e_1 -e_2\}\cup\{e_1 -e_i +e_{i+1},e_1 +e_i -e_{i+1}:
  2\leq i\leq k\}
\]
confined to $\mathbb{Z}^{k+1}_{\geq 0}$.
The number of walks of length $n$ from the
origin to the point $(n,0,...,0)$ staying within the nonnegative
octant equals the number of $n$-cell SYTs with at most $2k + 1$ rows.

If, additionally, the $e_1$ steps are confined to the hyperplane
spanned by the vectors $\{e_1, \dots, e_k\}$, then the number of paths
equals the number of $n$-cell standard Young tableaux with at most
$2k$ rows.

The number of $e_1$ steps is equidistributed with the odd column
statistic of standard Young tableaux. 
\end{theorem}

\begin{figure}\center
\includegraphics[height=4cm]{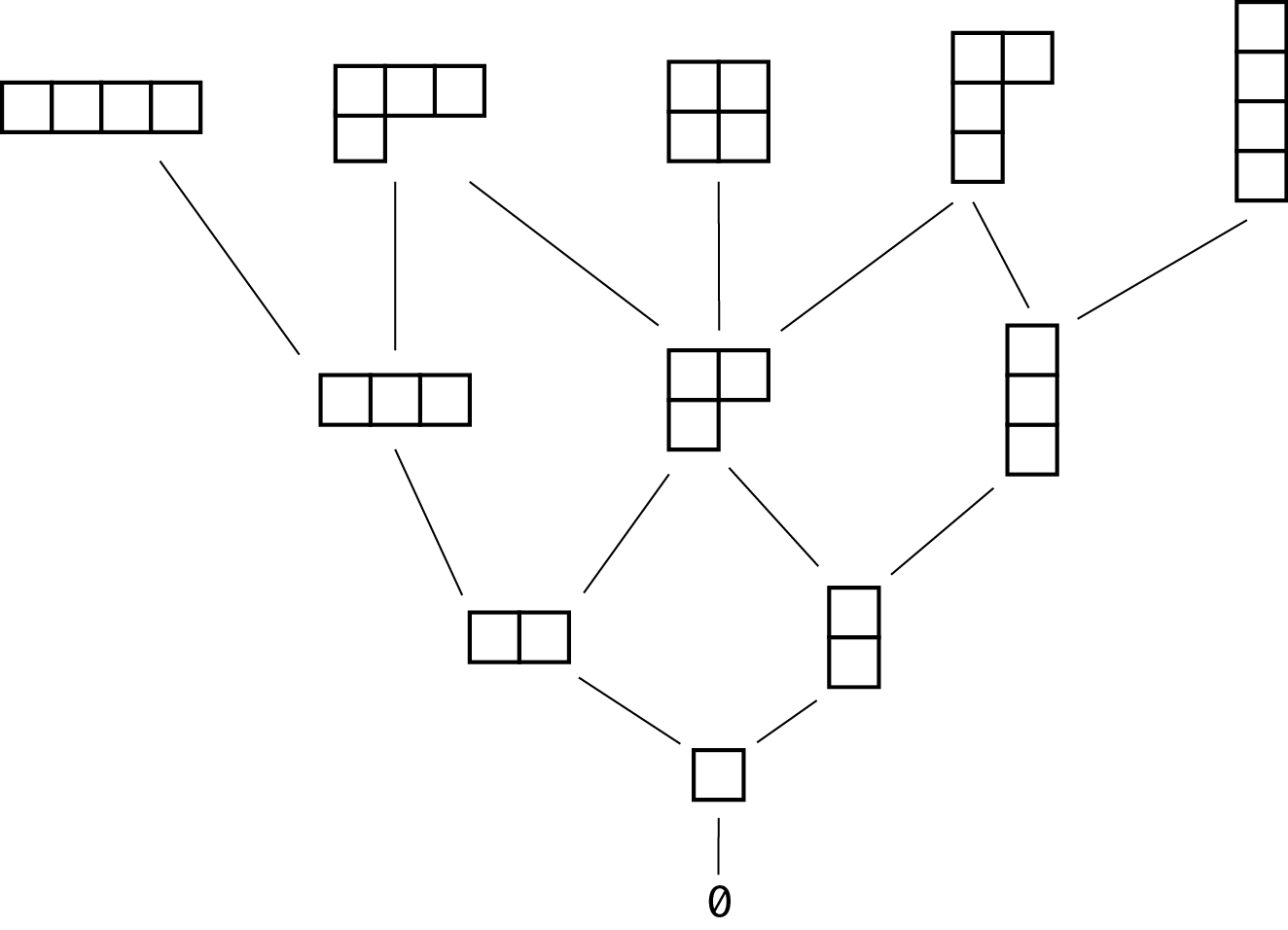}
\caption{The first few levels of Young's lattice of Ferrers diagrams}
\label{fig:YoungLattice}
\end{figure}

\subsection{Oscillating tableaux and arc diagrams}
\label{sec:axis} The set of Ferrers diagrams ordered by diagram
inclusion\footnote{Recall $\lambda\leq\mu$ means that
  $\lambda_i\leq\mu_i$ for all $i$} is called Young's lattice. Figure~\ref{fig:YoungLattice}
depicts the first few levels of its Hasse diagram.  We consider a sequence of
Ferrers diagrams as a walk on this
lattice. We consider three variants defined by restrictions on moving
up or down in the lattice (or not at all). The \emph{length\/} of a
sequence is the number of elements, minus one. (It is the number of
steps in the corresponding walk.)

An \emph{oscillating tableau} is simply a sequence of Ferrers diagrams
such that at every stage a box is either added or deleted. They were
popularized by their use in interpretations for representations of the
symplectic group~\cite{Sund90}.  If no diagram in the sequence is of
height $k+1$, we say that the tableau has its \emph{height bounded
  by~$k$}.  We recall two related families,
namely \emph{vacillating tableaux}, and \emph{hesitating
  tableaux}. The vacillating tableaux are \emph{even} length
sequences of Ferrers diagrams, written
$(\lambda^{(0)}, \dots, \lambda^{(2n)})$ where consecutive elements in
the sequence are either the same or differ by one square, under the
restriction that
$\lambda^{(2i)}\geq \lambda^{(2i+1)}$ and
$\lambda^{(2i+1)}\leq \lambda^{(2i+2)}$. The hesitating
  tableaux are even length sequences of Ferrers diagrams,
written $(\lambda^{(0)}, \dots, \lambda^{(2n)})$ where consecutive
differences of elements in the sequence are either the same or differ
by one square, under the following restrictions:
\begin{itemize}
\item if $\lambda^{(2i)} = \lambda^{(2i+1)}$, then $\lambda^{(2i+1)} < 
  \lambda^{(2i+2)}$ (do nothing; add a box)
\item  if $\lambda^{(2i)} > \lambda^{(2i+1)}$, then $\lambda^{(2i+1)} = 
  \lambda^{(2i+2)}$ (remove a box; do nothing)
\item if $\lambda^{(2i)} < \lambda^{(2i+1)}$, then $\lambda^{(2i+1)} >
  \lambda^{(2i+2)}$ (add a box; remove a box).
\end{itemize}  
Figure~\ref{fig:tableaux} shows examples from each of these classes.


\begin{figure}
\ytableausetup{smalltableaux}
\center
$${\tiny \left(\varnothing,  \varnothing,
\yt{\ },  \yt{\ }, \yt{\  \\ \ }, \yt{\   }, \yt{\ & \ },\yt{\ & \ },\yt{\ & \ }, \yt{\  & \  } , \yt{\  & \ \\ \ } \right)}$$
$${\tiny \left(\varnothing, \varnothing, 
\yt{\ },  \yt{\ }, \yt{\  & \ }, \yt{\ & \ \\ \ },\yt{\ \\ \ }, \yt{\  }, \yt{\  } \right)}$$
$${\tiny \left(\varnothing,  
\yt{\ },  \yt{\ \\ \ }, \yt{\  & \  \\ \ }, \yt{\  & \ }, \yt{\  & \
\\ \ },\yt{\ & \ \\ \ & \ }, \yt{\ & \ & \ \\ \ & \ },  \yt{\ & \ &
\   \\ \  }, \yt{\  & \ & \   },  \yt{\  & \ & \  \\ \ },\yt{\  & \
\\ \ } \right)}$$
 \caption{\textit{From top to bottom:} a vacillating tableau of length
   10; a hesitating tableau of length 8; an oscillating tableau of length
   11.  In each case, the height is bounded by $2$. From~\cite{Buetal16}}
\label{fig:tableaux}
\mbox{}\\
{\color{gray!50} \hrule}
\end{figure}

\ytableausetup{textmode,nosmalltableaux}

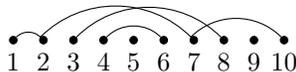
\begin{figure}
\centering
\begin{tikzpicture}[scale=0.4]
\foreach \x in {1,2,...,10}{
\node[pnt, label=below:{$\x$}](\x) at (\x,0){};}
\draw[bend left=45](1) to (2) to (7) to (10);
\draw[bend left=45](3) to (8);
\draw[bend left=45](4) to (6);
\end{tikzpicture}
\caption{An arc diagram representation of the set
  partition~$\pi=\{1,2,7,10\},\{3,8\},\{4,6\},\{5\},\{9\}$. It is both 3-noncrossing and 3-nonnesting. However, the
  subdiagram induced by $\{2,3,7,8,10\}$ is an enhanced
  3-crossing. Similarly, the subdiagram induced by $\{ 3,4,5,6,8\}$ is
  an enhanced 3-nesting.} 
\label{fig:expartition}
\end{figure}

Chen, Deng, Du, Stanley and Yan~\cite{Chetal07} described non-trivial
bijections between sequences of Ferrers diagrams and several combinatorial
families encoded in arc diagrams. \emph{Arc diagrams} are labelled
graphs under some degree and embedding restraints.  They can be used
to represent a variety of combinatorial classes, such as matchings and
partitions. Figure~\ref{fig:expartition} illustrates how to encode a
set partition as an arc diagram.

They generalize the Schensted algorithm, in some sense, by
describing how to parse arc diagrams, with each step defining a
tableau insertion or deletion. The result is bijection $\phi$ from
partitions to vacillating tableau. It is robust, and can be adapted to
other arc diagram classes. A key feature is sub-pattern avoidance
properties in the arc diagrams are mapped to height restrictions on
the tableaux. This echoes a key feature of the Schensted algorithm.

An arc diagram is noncrossing if no two arcs intersect. Noncrossing
set partitions are counted by Catalan numbers. In addition to
appearing in combinatorics, these diagrams arise in algebra, physics,
and free probability.  The notion of a crossing is generalized to a
\emph{$k$-crossing}, which denotes a set of $k$ arcs that each
mutually cross. Similarly, a \emph{$k$-nesting} refers to $k$ arcs
which mutually nest into a rainbow figure. More formally, let us
consider a set of~$k$ distinct arcs: $(i_1, j_1), \dots, (i_k,
j_k)$. They form a $k$-crossing if
$i_1<i_2<\dots<i_k<j_1<j_2<\dots<j_k$, and a $k$-nesting if
$i_1<i_2<\dots<i_k<j_k<\dots<j_2<j_1$. A slightly relaxed notion is
sometimes appropriate.  Enhanced $k$-nestings and $k$-crossings permit
the middle inequality to be an equality.  We are interested in the
classes of diagrams which \emph{avoid} such sub-diagrams, and consider
$k$-noncrossing diagrams (they contain no $k$-crossing) and
$k$-nonnesting diagrams (which contain no
$k$-nesting). Figure~\ref{fig:expartition} gives examples.

Theorem 3.2 of Chen \emph{et al.}~\cite{Chetal07} proves that given 
a partition $\pi$  of $\{1, \dots n\}$, and the vacillating tableau
$\phi(\pi)=(\emptyset=\lambda^{(0)}, \lambda^{(1)}, \dots,
\lambda^{(2n)}=\emptyset)$, the size of the largest crossing of $\phi$ is the
largest height of any $\lambda_i$, and the size of the largest nesting is the
maximum number of columns. 

Thus, by a simple transposition one deduces the corollary comprised of
the surprisingly nontrivial fact that $k$-noncrossing set partitions
of $\{1, 2, \dots, n\}$ are in bijection with $k$-nonnesting set
partitions of $\{1, 2, \dots, n\}$. We leave it to the reader to see
why a simple swap of sub-diagrams is not a bijection.

The bijection $\phi$ supports many generalizations. A diagram is said
to be \emph{open} if, in addition to closed edges, there are also half
edges. The notion of crossing extends naturally, given a fixed
convention of how the open edges lie.  The diagram $\mu$ in
Figure~\ref{fig:Collection} demonstrates an incomplete, or open,
matching. Remark there are not three edges that mutually cross, and
hence this diagram is also 3-noncrossing. An open diagram corresponds
to a sequence that does not necessarily end at an empty diagram--
rather it ends in a row shape, i.e. a partition with a single part:
$(m)$. The bijection $\phi$ is well defined in this case, although
Chen \emph{et al.} did not discuss how to interpret this, although
Burrill, Melczer and Mishna~\cite{BuMeMi15}. Xin~\cite{Xin10} consider
walks with any endpoint, and found generating functions for
palindromic sequences.

In a sequence of Ferrers diagrams of height at most $k$, each
diagram can be coded by a $k$-dimensional vector where each successive
entry is weakly decreasing. Thus, a sequence is an encoding of
a walk in $W_C(d)$. Thus we can view the bijection $\phi$ as a
nontrivial map between arc diagrams and walks. To end in a row shape
is the same as ending on the first axis.

This brings us back to our topic at hand, standard Young tableaux.
Which arc diagrams are associated to a standard Young tableau of
bounded height? The most plain interpretation is to view a standard
Young tableau as a sequence of Ferrers diagrams -- it starts at the empty diagram
and a single cell is added at each step, with the position indicated
by the entry in the tableaux. We could consider the arc diagram image of this
sequence, but we can do better.

Recall the Robinson-Schensted correspondence is a bijection between
pairs of standard Young tableaux and permutations. By restricting the
map to pairs of identical tableaux, it becomes a bijection between standard
Young tableaux and involutions. Involutions have a very natural arc
diagram representation! In this correspondence a tableau with a fixed
number of odd columns is mapped to an involution with precisely that
many fixed points. An involution is a partial matching, and they are
in bijection with oscillating tableaux. 

These pieces lead to the following result, a map between standard
Young tableaux and a class of walks that end on a boundary.
\begin{theorem}[Krattenthaler 2016~{\cite[Theorem 4]{Krat16}}; Burrill, Courtiel, Fusy, Melczer
  and Mishna 2016~{\cite[Theorem 1]{Buetal16}}; Okada 2016~{\cite[Theorem 1.2]{Okad16}}]
\label{thm:AxisWalks}
For $n,k\geq 1$, there is an explicit bijection between the standard
Young tableaux of size~$n$ with height at most $k$ and with $m$ odd
columns, and the simple walks of length $n$ staying in $W_C(k)$,
starting from the origin and ending at the point $m\,e_1$.
\end{theorem}

We have outlined one possible argument, but in fact three very
different proofs have been discovered. Above is the
argument of Burrill, Courtiel, Fusy, Melczer and
Mishna~\cite{Buetal16}. Krattenthaler~\cite{Krat16} determines a
different bijection using growth diagrams as an intermediary
object. It is direct, and gives several generalizations. It is self
contained in that it does not exploit the $\tau$
map. Okada~\cite{Okad16} uses a representation theoretic argument
starting from Pieri rules. He also has generalizations to
semi-standard tableaux. Gouyou-Beauchamps proved the $k=4$ case. 

One consequence of a lattice walk bijection is a generating
function expression.  Applying Grabiner and Magyar's formula gives the following
\[
Y_k(t)= \sum_{u=0}^{k-1}(-1)^u \sum_{\ell=u}^{2k-1-2u}(b_{\ell}) \det (b_{i-j}-b_{kd-i-j})_{0 \leq i \leq k-1, i\neq u, 1\leq j \leq k-1}.
\]
Remarkably, the infinite sum which arises from direct application of Grabiner and
Magyar's formula telescopes into a finite sum. We also use the
identity $b_{-k}=b_k$.

\subsection{Excursions in the Weyl chamber of type D}
Classically, there is a simple bijection between walks with steps from
$\{(1,1), (1, -1)\}$ which, in the first instance start and end on the
axis, with no further restriction (let us call them bridges), and in
the second instance start at the origin, and never go below the axis
and can end at any height (Dyck prefixes).  The first class
restricts free walks by restricting the possible end points, and the
second restricts free walks by restricting the size of the region. We
can view this as a tradeoff of restrictions. The simplest proof of
this result passes through a third object, marked Dyck paths, which
are objects in the intersection of both classes, but have additional
markings on some of the down steps which touch the axis. Consider the
following two maps, which both define bijections. A Dyck path with
marked down steps is mapped to a Dyck prefix, by flipping all of the
marked down steps into upsteps. The flipped marked steps become the
last step at that height in the walk, and the Dyck prefix ends at a
height equal to the number of marked steps.  On the other hand, we
could consider the entire segment of the walk ending at a marked down
step, starting at the nearest previous up step which touched the
axis. We can flip this entire segment across the axis. We do this for
every marked step to build a bridge. The marked intermediary
facilitates a straightforward proof, but how to mark in higher
dimensions?

The Weyl chamber of type $D$ is the following region:
\[W_D(k)\equiv\left\{ (x_1,x_2,\dots,x_k)\ \ :\ \ x_1 \geq x_2 \geq \dots
  \geq x_{k-1} \geq |x_k|\right\}.\] 
An \emph{axis-walk} is any walk starting
at the origin and ending on the $x_1$-axis. 
 \begin{theorem}[Courtiel, Fusy, Lepoutre and Mishna
   2018~{\cite[Theorem 20]{Coetal18}}]\label{theo:Wey}
   For~$k\geq 1$ and~$n\geq 0$, there is an explicit bijection between
   simple axis-walks of length $n$ staying in $W_C(k)$ and simple
  excursions of length $n$ staying in $W_D(k)$, starting from
   $(\frac 1 2,\dots,\frac 1 2,\frac 1 2)$, and ending at
   $(\frac 1 2,\dots,\frac 1 2,\frac {(-1)^n} 2)$. The ending
   $x_1$-coordinate of a walk from $W_C(k)$ corresponds to the number
   of steps that change the sign of $x_k$ in its bijective image.
\end{theorem}

The proof of this result is a generalization of the marked Dyck path
example. The analog to the Dyck prefix is the axis walk in the orthant
space $W_C(k)$, and the analog to the bridge is the excursion in the
larger space $W_D(k)$. The intermediary class, marked excursions in
$W_C(k)$, are defined by a less straightforward process. For the first
bijection, Courtiel~\emph{et al.} map the axis walks to open arc
diagrams. The open arcs are removed, but their position is marked. The
inverse bijection is applied and a marked excursion results. The
second bijection is slightly more complicated, and the reader is
referred to the article for details.

A corollary of this result is another lattice model in
bijection with standard Young tableaux of bounded height.

\begin{corollary}[Courtiel, Fusy, Lepoutre, Mishna 2018~{\cite[Corollary 21]{Coetal18}}]
\label{thm:Excursions}
For $n,k\geq 1$, there is an explicit bijection between the standard
Young tableaux of size $n$ with height at most $k$, and the simple
walks of length $n$ staying in $W_D(k)$, starting from $(\frac 1
2,\dots,\frac 1 2,\frac 1 2)$, and ending at $(\frac 1 2,\dots,\frac 1
2,\frac {(-1)^n} 2)$. The number of odd columns corresponds to the
number of steps that change the sign of $x_k$.
\end{corollary}

It seems very likely that an explicit bijection could be found between
these, and the lazy walks. This would answer the question of
Zeilberger, and likely reveal a model for the odd cases. 

There is also potential for another generating function formula by
applying the results of  Grabiner and Magyar.

\subsection{A collection of bijections}
\label{sec:Collection}
We summarize some combinatorial classes in bijection with standard
Young tableaux of bounded height.  Figure~\ref{fig:Collection} illustrates some of these
classes.

\begin{theorem}
The set of standard Young tableaux of size~$n$ with height bounded
by~$k$ and~$m$ odd columns is in bijection with each of the following sets: 
\begin{enumerate}
\item\label{invo} The set of involutions of size $n$ with $m$ fixed points and no decreasing subsequence of length $k + 1$;
\item\label{osc} The set of oscillating tableaux of size~$n$ with height bounded
by~$k$, which start at the empty partition and end in a row shape
$\lambda=(m)$;
\item\label{matching} (if $k$ is even) The set of open
matching diagrams of length $n$,  with $m$ open arcs and with no
$(k/2+1)$-crossing;
\end{enumerate}
\end{theorem}

Gil, McNamara, Tirrell and Weiner~\cite[Theorem 1.1]{Gietal17} proved
the equivalence between Class~\eqref{matching} in the theorem, and a
set of a Dyck paths where each groups of ascent steps are decorated by
a connected matching of a related size.  A main feature of their
construction is the following identity.  Suppose $P_k(t)$ denotes the
generating function for the number of $k$-noncrossing perfect
matchings on $[2n]$, then
\[Y_{2k-1}(t) = \frac{1 + P_k(t^2(1-t)^{-2})}{1-t}.\]

\section{Generating function expressions}
\label{sec:diag}
The generating functions of lattice walks in cones with symmetries in
the stepset can often be expressed as a subseries extraction of a
Laurent expansion of a rational function. This generalizes the result
on reflectable walks of Gessel and Zeilberger, and is at the heart of
the orbit sum method of Bousquet-M\'elou and Mishna~\cite{BoMi10}
which handles a wide class of lattice walks restricted to the quarter
plane. In their work, they start with a natural combinatorial
recurrence on lattice walks, which translates to a functional equation
satisfied by the generating function. For some models it is sufficient
to take a weighted sum (the namesake orbit sum) of the equations, and
isolate the target generating function with a sub-series extraction.

In this context, the extraction operates on iterated Laurent
series.  For a function
$F(x_1, \dots, x_d;t)\in \mathbb{C}(x_1, \dots, x_d)(t)$ which is
analytic at the origin, we denote by $[x_1^{k_1} x_2^{k_2}\dots
x_d^{k_d}t^n] F(x_1, \dots, x_d;t)$ to be the coefficient of the term $x_1^{k_1} x_2^{k_2}\dots
x_d^{k_d}t^n$ in the Laurent expansion. In this context we view the
objects to be series in $t$ variable. 
We recall the terminology that
several authors use for the special case of the constant term with
respect to a set of variables:
$\operatorname{CT}_{x_1\dots x_d}F (x_1, \dots, x_d; t)$ is the
coefficient of the term $x_1^0\dots x_d^0$ in the Laurent series
expansion of $F$ at the origin. It is a series in the remaining
variables. It can be obtained by incrementally determining the
constant term with respect to a single variable. Here we shall take
the order of the variables as they are listed, for example.  Such a coefficient
extraction can be written as a Cauchy integral.

A related operator is the \emph{diagonal}. The (central) diagonal
$\Delta F(x_0, x_1, \dots, x_d)$ of a formal power series is the
univariate subseries defined
\[
  \Delta F(x_0, \dots, x_d)=%
  \Delta \sum_{i_0, \dots, i_d} f(i_0,\dots, i_d) x_0^{i_0}\dots z_d^{i_d} %
  \equiv \sum_n f(n,\dots, n)\,t^n.
\]
Lipshitz~\cite{Lips88} proved that the diagonal of a D-finite
function\footnote{A (multivariate) function is D-finite if the set of
  all its partial derivatives spans a vector space of finite
  dimension.}  is also D-finite. Since rational functions are
D-finite, diagonals of rational functions are also D-finite. The
expressions we consider here are all diagonals of rational functions,
and hence are all D-finite by construction. This provides an alternate
proof of the D-finiteness of $Y_k(t)$. 

Bousquet-M\'elou~\cite{Bous11} used an orbit sum to directly derive
generating function expressions for standard Young tableaux of bounded
height.  Her generating functions also mark the number of odd columns,
and other parameters are also readily accessible.  In that same paper,
she determines a nice functional equation proof of the hook length formula.
The starting point is a simple recursive construction for standard
Young tableaux: a tableau of size $n+1$ is obtained from a tableau of
size $n$ by adding a cell to the $j$-th row, unless the $j$-th row is
already the same length as the $j-1$ row. This translates into 
a very straightforward functional equation for the generating function
\[
F(u)\equiv \sum_{\substack{\lambda\vdash n\\
  \ell(\lambda)\leq k}} f^\lambda u_1^{\lambda_1}u_2^{\lambda_2}\cdots.
\]
The equation is defined in term of $F_j(u)$, the generating function for those
standard Young tableaux such that $\lambda_{j-1}=\lambda_j$:
\[
F(u)= 1+u_1F(u) +\sum_{j=2}^k u_j(F(u)-F_j(u)).
\]
She applies a kernel method argument to this functional equation to
recover MacMahon's formula for $f^\lambda$.

She has generating function expressions in her Propositions~9
and~10. It is her Proposition~11 that accounts for the additional
parameter of interest. We note that she works with the ordinary
generating function, 
\[
\widetilde{Y}_k(u;t) \equiv \sum_n \sum_{\substack{\lambda\vdash n\\
  \ell(\lambda)\leq k}} f^\lambda\, u^{\# \text{odd columns in }\lambda}\,t^n.
\] 
\begin{theorem}[Bousquet-M\'elou 2011~\cite{Bous11}]
For $k=2\ell$, the ordinary generating function of standard Young
tableaux with height bounded by~$k$ where~$t$ marks the length and
$x_1$ marks the number of odd columns is 
\begin{equation}\label{eqn:MBM}
\widetilde{Y}_k(x_1;t) = \operatorname{CT}_{x_2,\dots, x_k} \frac{-\det(x_j^{-i}-x_j^i)_{1\leq i, j\leq
    \ell}}{1-t(x_1+x_1^{-1}+\dots+x_\ell+x_\ell^{-1})}\frac{1}{x_2x_3^2\dots
x_\ell^{\ell-1}}.
\end{equation}
\end{theorem}
The orbit sum strategy can be used to determine an expression for the
ordinary generating function for simple walks in $\mathcal{W}_C(k)$
ending on the axis. This gives the following result. 
\begin{theorem}[Burrill, Melczer and Mishna 2015~\cite{BuMeMi15}]
\label{thm:osc}
The ordinary generating function, $\widetilde{Y}_k(t)$ for Young tableaux of
height bounded by~$k$ satisfies the formula
 formula
{\small
\begin{align*}
   &t^{2k-1}\widetilde{Y}_{2k}(t) =\\
& -\Delta \left[ \frac{z_0^{2k-1}(z_3 z_4^2 \cdots
    z_k^{k-2})(z_1+1)\prod_{1\leq j<i \leq k} (z_i-z_j)(z_iz_j-1)
    \cdot \prod_{2 \leq i \leq k} (z_i^2 -1)}{1-z_0(z_1\cdots z_k)(z_1 +
    z_1^{-1} + \cdots z_k + z_k^{-1})}  \right]. 
\end{align*}}
\end{theorem}


These formulas are not necessarily easier to compute or interpret than
the others, but they do have a few implications. They do provide a
proof that the generating function is D-finite. Furthermore, we see a
second proof that $y_k(n)$ can be expressed as binomial sums, and
potentially a different route to obtain such expressions, using the
work of Bostan, Lairez and Salvy~\cite{BoLaSa17}. They may be
suitable for analysis by the methods of Pemantle and
Wilson~\cite{PeWi13}, which determine asymptotic formulas for
coefficients of functions which are expressed as diagonals of rational
functions.

Since it is the generating function for Motzkin numbers,
$\widetilde{Y}_3(t)$ is algebraic. Because the asymptotics of the
coefficients are incompatible with algebraicity, $\widetilde{Y}_4(t)$
is not. For which $k$ is the ordinary generating function $\widetilde{Y}_k(t)$
algebraic? 

\subsection{Differential equations}
\label{sec:bounds-on-DE}
Once we know that a generating function is D-finite, it
is natural to ask about differential equations that it
satisfies. Bergeron, Favreau and Krob~\cite{BeFaKr95} generated many
conjectures about the order of the differential equation satisfied by
$Y_k(t)$ from computer experiments, and some analysis of its
expression as a determinant of a matrix of modified Bessel functions.

Proposition 1 in~\cite{BeGa00} states the
dimension of the vector space over the field $\mathbb{C}(t)$ of rational
functions in $t$ spanned by $Y_k(t)$ and all its derivative is bounded by
$\lfloor\frac{k}{2}\rfloor$. 

\begin{conjecture}[Bergeron, Favreau, Krob 1995~\cite{BeFaKr95}] 
For each $k$, there are polynomials $p_m(t)$ of
  degree at most $\lfloor\frac{k}{2}\rfloor$ such that $Y_k(t)$ is a
  solution of a linear differential equation order at most
  $\lfloor\frac{k}{2}\rfloor+1$ with coefficients $p_m(t)$. 
\end{conjecture}
They have verified this conjecture for $k\leq 11$. 

We can use recent advances in symbolic computation and the diagonal
expression in Theorem~\ref{thm:osc} to bound the order and degree of
the differential equations satisfied by the ordinary generating
function $\widetilde{Y}_k(t)$. In particular, the following theorem of
Bostan, Lairez, and Salvy~\cite{BoLaSa13} can be explicitly applied.
\begin{theorem}[Bostan, Lairez, and Salvy~\cite{BoLaSa13}]
  Let
  $R(z_1,\dots,z_d,t) = A(\bz,t)/B(\bz,t) \in
  \mathbb{Q}(t)(z_1,\dots,z_d)$, be a rational function with
  multidegree bounds
\[ n_\bz := \max(\deg_\bz B, \deg_\bz A + d + 1) \qquad\qquad n_t := \max(\deg_t A, \deg_t B). \]
Then there exists an annihilating differential equation $\mathcal{L}$ for the integral 
\[P(t) := \oint_\gamma R(\bz,t) d\bz,  \]
where $\gamma$ is any $n$-cycle in $\mathbb{C}^n$ on which $R$ is continuous when $t$ ranges over some connected open set $U \subset \mathbb{C}$ (note that $\mathcal{L}$ is independent of $\gamma$).  Furthermore
the order of $\mathcal{L}$ is at most $n_\bz^d$ and 
the degree of $\mathcal{L}$ is at most $\left(\frac{5}{8}n_\bz^{3d} + n_\bz^d\right)e^d n_t$.
\end{theorem}
We write a diagonal as an integral:
\begin{equation}
\Delta F(\bz,t)=\left(\frac{1}{2\pi i}\right)^d \int_\gamma \frac{F(z_1,z_2/z_1,z_3/z_2,\dots,z_d/z_{d-1},t/z_d)}{z_1 z_2 \cdots z_d} d\bz
\end{equation}
by the multivariate Cauchy residue theorem, for an appropriate $n$-cycle $\gamma$ around the origin.

In the case of $\widetilde{Y}_k(t)$, we compute an upper bound of $2d^2-3d+1$ on
the total degree of the denominator in the $\bz$ variables and a bound
of $2d^2-4d-2$ on the degree of the numerator and implies the
following result, computed by Melczer.
\begin{proposition}
  The generating function for the number of standard Young tableaux of
  height at most $k$ satisfies a linear differential equation of order
  at most $(2k^2-3k+1)^k$ and degree at most 
\[\left(\frac{5}{8}(2d^2-3d+1)^{3d}+(2d^2-3d+1)^d\right)e^d.\]
\end{proposition}

\subsection{Asymptotics}
Regev~\cite{Rege81} determined asymptotic expansions for the number of
Young tableaux of bounded height. He explicitly deduced asymptotics
for the $k=3$ case as a particular extraction:
\[(n+1)y_3(n)=[x^{-1}]  (x+1/x+1)^{n+1}\sim \sqrt{\frac{3}{8}\pi}\cdot
  \frac{1}{\sqrt{n}} 3^n.
\]

By a slightly more general argument, he showed 
\begin{equation}\label{eq:asympt}
y_{2k}(n)\sim_{n\to\infty}
(2/\pi)^{k/2}(2k)^n(k/n)^{k(k-1/2)}\prod_{i=0}^{k-1}(2i)!.
\end{equation} 

The asymptotics of lattice walks in cones have been well studied.  The
tour de force of Denisov and Wachtel~\cite{DeWa15} describes a
collection of very comprehensive results. The formulas given
in~\cite[Theorem 6]{DeWa15} should be applied here, for example to the
lazy walks and to the excursions in $W_D(k)$.

Grabiner~\cite{Grab06} used his formulas for walks in Weyl chambers to
find the asymptotics of the probability that a randomly chosen
standard Young tableau of size~$n$ with at most~$t$ rows contains a
given subtableau. This is equivalent to counting walks that have
visited a particular point-- there might be similar results to be
extracted by considering other lattice models. 

The asymptotic formula that have been described so far are valid for
fixed $k$, as $n$ tends to infinity. It is open to develop formulas
when $k$ is a function of $n$.

\section{Restricting increasing subsequences in permutations}
\label{sec:permutations}
Motivated by algebraic interpretations, Regev~\cite{Rege81} considered the quantity
\[
y_{k}^{(\beta)}(n) \equiv \sum_{\lambda\in\mathcal{P}_k} (f^\lambda)^\beta
\]
What can be said of the generating functions
\[
\sum y_{k}^{(\beta)}(n) t^n\quad  \sum y_{k}^{\beta}(n) t^nu^\beta?
\]

When $\beta=2$, this counts permutations with restricted longest
increasing subsequence, and it is very well studied, with many
relevant connections to lattice walks~\cite{Wilf95, Wilf92}. Define
$u_k(n)$ to be the number of pairs of Young tableaux of the same shape
with at most $k$ rows. By the Robinson-Schensted correspondence, this
is the number of permutations in $\mathfrak{S}_n$ with no
$(k+1)$-increasing subsequence. For every $k \geq 1$,
Gessel~\cite{Gess90} proved the formula
\begin{equation}
\sum_{n\geq 0}\frac{u_k(n)}{n!^2t^{2n}}=\det\left(b_{i-j}\right)_{1\leq i,j\leq k},
\end{equation}
where the $b_j$ are the Bessel function evaluations defined earlier. 

A combinatorial proof of this expression has been given by Gessel
\emph{et al.}~\cite{GeWeWi98} via simple walks ending at so-called
Toeplitz points, and Xin~\cite{Xin10} did a combinatorial derivation
based on arc diagrams, and yet another constant term extraction. The
permutations satisfy a nice combinatorial recursion, and, in a manner
similar to her solution for involutions,
Bousquet-M\'elou~\cite{Bous11} determines functional equations that
can be resolved using a kernel method approach.

She determined~\cite[Proposition 13]{Bous11} the following expression for the ordinary generating
function of permutations avoiding the pattern~$1\, 2 \dots m \, (m+1)$:
\begin{equation}
U_{k}(t)\equiv\sum_{n\geq 0}u_k(n)t^{2n} = \operatorname{CT}_{x_1, \dots, x_m}
\frac{\det((x_j-x_{j-1})^{i-j})_{1\leq i,j\leq m}}{1-t\left(\sum
    \frac{1}{x_j-x_{j-1}}\right)} \cdot \sum_{i=0}^m \prod_{j=1}^i \frac{x_j}{1-x_j}.
\end{equation}

Bergeron and Gascon found the differential equations satisfied by the
exponential generating functions for $k<11$. It remains open to answer
if $(y_{k}^{(\beta)}(n))$ is the counting sequence for any easily
characterizable combinatorial family, perhaps as a restricted or
decorated family of permutations.  The argument of
Gessel~\cite{Gess90} applies to prove that the sequence is P-recursive
for positive integer $\beta$, and probably also a diagonal of a
rational function.  Can we usefully bound the annihilating
differential operators or determine asymptotic formulas for arbitrary
$\beta$?

Wilf~\cite{Wilf92} deduced $U_{2k}(t)=Y_{2k}(t)Y_{2k}(-t)$. From this,
it follows
\[
\binom{2n}{n} u_{2k}(n)=\sum_r\binom{2n}{r}(-1)^ry_{2k}(r)y_{2k}(2n-r).
\]
Is there a lattice path interpretation of this identity? Are there
identities for other $\beta$ values?

\section{Other Directions}
\label{sec:Remarks}
\subsection{Using Kronecker coefficients}
The Kronecker product of symmetric functions gives an important
connection to representation theory. In particular, the
Kronecker product of two symmetric functions in the Schur function basis
determines the multiplicities of irreducible characters in this tensor product.

The following Schur function identity for the Kronecker product of two
Schur functions (denoted by $*$) was shown by Brown, van Willigenburg and Zabrocki~\cite{BoVaZa10}: 
\[
s_{(n, n-1)}*s_{(n, n-1)}=\sum_{\substack{\lambda\vdash 2n-1\\ \ell(\lambda)\leq 4}}
s_\lambda. 
\]
Tewari~\cite{Tewa15} manipulates this formula  to deduce a closed form for the number of Young tableaux with height exactly 5, under the additional
constraint of $\lambda_5=1$. His Theorem 7.4 is a simple sum of
Motzkin numbers (recall $M_n=y_3(n)$) and pairs of Catalan numbers
(recall also the formula of $y_4$):
\begin{equation*}
\label{eqn:Tewari}
\sum_{\substack{\lambda\vdash 2n\\ \ell(\lambda)=5 \\\lambda_5=1}}
f^\lambda =\left(\frac{n(n+2)}{2n+1} \right) C_n C_{n+1}-C_{n+1}^2+M_{2n},
\end{equation*}
\begin{equation*}
\label{eqn:Tewari2}
\sum_{\substack{\lambda\vdash 2n-1\\ \ell(\lambda)= 5 \\\lambda_5=1}}
f^\lambda =\left(\frac{(n+1)}{2} \right) C_n^2-C_nC_{n+1}+M_{2n-1}.
\end{equation*}
It could be straightforward to find a combinatorial interpretation of his formula,
which then could be generalized to higher dimensions. It also suggests
that perhaps some of the larger summations could be simplified. This
could be key to finding a useful expressions for $y_6(n)$.

\subsection{Other classes in bijection}
As Gouyou-Beauchamps noted in~\cite{Gouy89}, the numbers that appear
in $y_4(n)$ also appear in the enumeration of planar maps and
alternating Baxter permutations.  Baxter permutations are a class of
pattern avoiding permutations that are very combinatorially rich. They
also have a bijection to lattice path models, using the machinery of
arc diagrams as an intermediary class. It might be possible to
directly connect these classes. 

\subsection{Shadow diagrams}
Some of the lattice walk bijections use tableau insertion and deletion
in more than one stage. Perhaps there exist more economical, or direct
bijections. The shadow diagrams of Viennot may play an important role
in such a simplification.

\subsection{Random tableaux}
There have been several recent works on generating random walks. One
application is to convert this to a generator for random tableaux.
Which of the above bijections has smallest complexity?

Grabiner~\cite{Grab06} was able to use lattice walk results to
determine distributions of subtableaux in Young tableaux. It seems
that there should be more results along this vein with each of the
lattice walk representations. 

\subsection{Semi-standard Young tableaux}
Okada~\cite{Okad16} proved some results connecting
counting sequences of generalized oscillating tableaux and
semi-standard tableaux using techniques from representation theory. In
particular, his Theorem 5.3 is a list of results on equinumerous
classes of tableaux of bounded height that are well suited for more
bijective explanations. 


Krattenthaler's results~\cite{Krat16} on semi-standard tableaux
replace single steps in oscillating tableaux with jumps by horizontal
strips. Perhaps they can connect here using lattice walks with longer
steps or diagonal steps. In any case, these two works should be connected more
explicitly.

\section*{Acknowledgements}
The author is grateful to MSRI for travel support to participate in
the 2017 AWM session. This expository work was inspired by that
meeting. I am grateful for the patience and wisdom of the anonymous
referees. The author's research is also partially funded by NSERC
Discovery Grant RGPIN-04157.



\end{document}